\documentclass[12pt]{article}       
\usepackage{a4}
\usepackage{amsmath,amssymb,amsthm} 
\newcommand{\N}{\mathbb N}        
\newcommand{\Z}{\mathbb Z}        
\newcommand{\Q}{\mathbb Q}        

\theoremstyle{plain}               
\newtheorem{definition}{Definition}[section]

\newtheorem{proposition}[definition]{Proposition}
\newtheorem{corollary}[definition]{Corollary}

\theoremstyle{definition}          

\newtheorem{example}[definition]{Example}

\begin{document}

\title{Mersenne Binomials and the Coefficients of the
non-associative Exponential}

\author{by L.~Gerritzen}
\date{}
\maketitle
\begin{abstract}

The non-associative exponential series $exp(x)$ is a power series
with monomials from the magma $M$ of finite, planar rooted trees.
The coefficient $a(t)$ of $exp(x)$ relative to a tree $t$ of
degree $n$ is a rational number and it is shown that

$$\hat{a}(t) := \frac{a(t)}{2^{n-1}\cdot \prod^{n-1}_{i=1}(2^i -
1)}$$

is an integer which is a product of Mersenne binomials. One
obtains summation formulas

$$\sum \hat{a}(t) = \omega(n)$$

where the sum is extended over all trees $t$ in $M$ of degree $n$
and $$\omega(n) = \frac{2^{n-1}}{n!} \prod^{n-1}_{i=1} (2^i -
1).$$

The prime factorization of $\omega(n)$ is described. The sequence
$(\omega(n))_{n \ge 1}$ seems to be of interest.
\end{abstract}

\bigskip\bigskip

{\bf{\large{Introduction}}}

\bigskip
  For a natural number $n,$ we denote by
$M_n$ the n-th Mersenne number $2^n-1.$

The Mersenne factorial $n!_M$ is defined to be the product
$\displaystyle \prod^n_{i=1} M_i$ of all Mersenne number
$M_1,...,M_n$ while the Mersenne binomial ${n \choose r}_M$ is
defined to be $${n!_M}\over {r!_M(n-r)!_M}$$ for natural numbers
$r$ between $0$ and $n.$

\medskip
In quantum calculus, see for instance [KC], there are the notions
$[n]$ for the $q$ - analoque of $n \in \N,[n]!$ for the $q$ - factorial
and $\big[{n \atop r}\big] $ for the $q$ - binomial coefficients which
are polynomials in the
variable $q$. If $f(2)$ denotes the value of a polynomial $f$ after substituting
2 for $q$, then $[n](2) = M_n, [n]!(2) = n!_M$ and $\big[{n \atop r}\big](2)
 = {n \choose r}_M.$

\bigskip
Let $M$ be the magma with unit $1_M$ freely generated by a single
element $x$. It can be identified with the set of all finite planar
binary rooted trees, see [DG].

\medskip

Let $\Q\{\{x\}\}$ the $\Q$- algebra of power series with
monomials from $M$. It was proved in [DG] that there is a unique
series $exp(x)\in\Q\{\{x\}\}$ such that $exp(x) = 1 + x + $higher terms
and $$exp(x) \cdot exp(x) = exp(2x)$$ see also Proposition (3.2).
Let $a(t)$ be the coefficient of the non-associative exponential
$exp$ relative to the planar binary rooted tree $t$ and put
$$\hat{a}(t) = n! \cdot \omega(n)\cdot a(t)$$ where $$\omega(n):=
{{2^{n-1}(n-1)!_M}\over {n!}}.$$

\medskip
A main result of this article states that all $\hat{a}(t)$ are
integers which are obtained as products of Mersenne binomials. The
projection of the non-associative exponential to the classical one
leads to the decomposition $$\sum \hat{a}(t) = \omega(n)$$ if the
sum is extended over all trees of degree $n.$

\medskip
In section 1 the notions of Mersenne order and Wieferich exponents
at odd primes is defined. One gets results about the prime
factorization of Mersenne numbers.

\medskip

In section 2 we show that all $\omega(n)$ are integers by
computing the order $ord_p \omega(n)$ for all primes $p.$


The applications about the coefficients $a(t)$ of $exp(x)$ are
deduced in section 3.

\medskip

I would like to thank Doron Zeilberger for helpful hints.

\begin{section}{Mersenne orders}

We call $M_n = 2^n-1$ the n-th Mersenne number for $n \in \N, n
\ge 1.$

This is in contrast to the definition by D. Shanks, see [S], Chap.
1, which is also widely used in which $M_n$ is called Mersenne
number only if $n$ is a prime.

Let $p$ be an odd prime.

\begin{definition}
The Mersenne order $v(p)$ at $p$ is the smallest number $r \in \N,
r \ge 1,$ such that $p$ divides $M_r.$

\end{definition}

 As $2^{v(p)} \equiv 1 \ mod \ p$ and $v(p)$ is the order of
the class of $2$ in $(\Z/p\Z)^*,$ we get that $v(p)$ divides the
order of the multiplicative group $(\Z/p\Z)^*$ which is $(p-1).$

Also $p$ divides $M_n,$ if and only if $n$ is a multiple of
$v(p).$

Also $p < 2^{v(p)}$ from which follows that $log_2(p) < v(p).$

It has the consequence that $\displaystyle lim_{p \to \infty}
v(p)= \infty.$

\begin{definition}
Let $\varepsilon(p) := ord_p(M_{v(p)})$ for any odd prime; it is
called the Wieferich exponent at $p.$
\end{definition}

 \medskip

 Always $\varepsilon(p) \ge 1$ and $p$ is called
a Wieferich prime if $\varepsilon (p) \ge 2.$

Up to now only two Wieferich primes are known, namely 1093 and
3511 and it was checked that these are the only Wieferich primes $
< 4.10^{12}$ in [CDP].

The Mersenne order at 1093 is 364 = ${{1092} \over {3}} = 2^2
\cdot 7 \cdot 13$ and the Mersenne order at 3511 is 1755 =
${{3510}\over {2}} = 3^3 \cdot 5 \cdot 13.$

The Wieferich exponent for both primes is equal to 2.

\begin{proposition}

Let $p$ be an odd prime, $n \in \N, n \ge 1.$

\begin{itemize}
\item[(i)] $ord_p(M_n) = 0,$ if $n$ is not a multiple of $v(p).$
\item[(ii)] If $n = v(p)\cdot p^r \cdot k$ and $p$ does not divide
 $k$, the  $ord_p(n) = r$ and $ord_p(M_n) = \varepsilon (p) + r$
 \end{itemize}
\end{proposition}

\begin{proof}
Statement (i) is well-known.

Let now $n = v(p) \cdot p^r \cdot k$ as in (ii).

Let $m$ be an integer and consider the group $(\Z/p^m\Z)^*$ of
units of the ring $\Z/p^m\Z.$ Let $U_m$ be the subgroup of
$(\Z/p^m\Z)^*$ of class congruent to 1 modulo $p.$

Then $U_m$ has order $p^{m-1}$ because any element $\alpha$ in
$U_m$ has a unique representation by a number $a = 1 +
\displaystyle \sum ^{m-1}_{i=1} a_i p_i$ with $a_i \in
\{0,...,p-1\}.$ The order of $\alpha$ is $m - min(\{i: a_i \ne
0\}, m)$ because if $a = 1 + p^\varepsilon \cdot x, x $ any number
not divisible by $p$, then $$a^p = 1 + p^{\varepsilon + 1}\cdot
y$$ where $$ y = x + \displaystyle \sum^p_{i=2}(^p_i) x^i p ^{i
\varepsilon - \varepsilon - 1}$$ which shows that $y \equiv x \
mod \ p$ as $p > 2.$

Let now $\alpha = \bar {2}^{v(p)}$ where $\bar {2}$ is the class
of 2 in $\Z/p^m\Z.$

Then $\alpha \in U_m$ and $\alpha$ is represented by a number $a =
1 + p^{\varepsilon(p)} \cdot x,$ where $x$ is not divisible by
$p.$ Then  $\alpha^{p^r}$ is represented by a number $$1 +
p^{\varepsilon(p) + r} \cdot y$$ with $y$ not divisible by $p.$

It follows that $ord_p M_n = \varepsilon(p) + r$ because the order
of $\alpha^{p^r}$ is equal to the order of $\alpha^{p^rk}.$
\end{proof}

\begin{corollary}
$$M_n = \prod_n p^{\varepsilon(p) + ord_p(n)}$$

where the product is extended over all odd primes $p$ for which
the Mersenne order $v(p$) divides $n.$
\end{corollary}

\end{section}
\begin{section}
{Factorial Mersenne quotients}

$F(n) := \prod^n_{i=1} M_i$ is called the Mersenne factorial of $n
\in \N;$ it is also denoted by $n!_M.$

\medskip

Let $$\omega(n) := {{2^{n-1} \cdot F(n-1)}\over {n!}}$$  for
$n \ge 1;$ it is called the factorial Mersenne quotient at $n.$
\medskip
 Let $m \in \N, m \ge 1,$ and $p$ be a prime. There is a
unique p-adic expansion $$m = \sum^r_{i=0}m_i p^i$$ with $m_i \in
\{0,1,...,p-1\}.$ Let $$d_p(m) := \sum^r_{i=0}m_i.$$

\begin{proposition}
\begin{itemize}
\item[(i)] $ord_2\omega(n) = d_2(n) - 1$
\item[(ii)]If p is an odd prime, then
$$ord_p(\omega(n)) = \varepsilon(p)\cdot m - {{(n-d_p(n)) - (m-d_p(m))}\over
{(p-1)}}$$
where $v(p)$ is the
Mersenne order at $p$,
$m$ is the largest integer $ \le (n-1)/v(p)$
\end{itemize}
 and
$\varepsilon(p)$ is the Wieferich exponent at $p$.
\end{proposition}

\begin{proof}
\begin{itemize}
\item[1)]  $m
\equiv d_p(m) \ mod \ (p-1)$ because $m \equiv \displaystyle \sum
^r_{i=0} m_i \ mod \ (p-1).$ It is well known that $$ord_p(n!) =
{{n-d_p(n)} \over {(p-1)}} = \sum ^{< \infty}_{i=1} \ [{{n} \over
{p^i}}]$$
\end{itemize}
\begin{itemize}
\item[1)] $ord_2(\omega(n)) = (n -1) - ord_2(n!)$ because $F(n)$ is
odd. Thus $$ord_2(\omega(n)) = d_2(n) - 1 \ge 0$$ for $n \ge
1.$
\item[2)] Let now $p$ be an odd prime. Then $$ord_pF(n-1) =
\sum^{n-1}_{i=1} ord_p M_i.$$ If $i= k \cdot p^r v(p)$ and $k$ is
not divisible by $p,$ then $$ord_pM_i = \varepsilon(p) + r)$$ and
$ord_p M_i = 0$ if $p$ is not divisible by $p.$ $$ord_p(F(n-1)) =
\sum^m_{j=1} ord_p M_{jv(p)}$$ where $m_p(n) =
m:=[{{n-1}\over {v(p)}}]$ is the largest integer $\le (n-1)/v(p)$.
It follows from Proposition 1.3 that $$ord_p
M_{j\sigma(p)} = \varepsilon(p) + ord_p(j)$$ as $ord_p(j) =
ord_p(jv(p)).$ As $$\displaystyle \sum^m_{j=1} ord_p(j) = ord_p(m!)
= {{m-d_p(m)}\over {p-1}}$$ we get the formulas $$ord_p F(n-1) =
\varepsilon(p) \cdot m + {{m-d_p(m)}\over {p-1}}$$and
$$ord_p(\omega(n) = \varepsilon(p)m - {{(n-d_p(n)) - (m-d_p
(m))}\over {p-1}}$$
 \end{itemize}
 \end{proof}

 \begin{proposition}
$\omega(n) \in \N.$
\end{proposition}

\begin{proof}
We will show now that $ord_p \omega(n) \ge 0.$ From the
formula in 3) we get $$ord_p(\omega(n)) \ge m - {{(n-d_p(n)) -
(m-d_p(m))}\over{p-1}}$$ $$\ge m - {{n-d_p(n)}\over {(p-1)}}\ge
[{{n-1}\over {p-1}}] - {{n-d_p(n)} \over {p-1}} \ge 0$$ because $m
- d_p(m) \ge 0, m \ge {{n-1} \over {p-1}},$ $n - d_p(n) \le n-1$
and ${{n-d_p(n)} \over {p-1}} \in \N.$
\item[(3)] From 1) and 2) we get that $ord_p(\omega(n-1)) \ge 0$ for
all primes $p.$ As $\omega(n-1)$ is a rational number, it follows
that $\omega(n-1)$ is an integer.

\end{proof}

\begin{example}
\begin{itemize}
\item[(1)] $ord_7 \omega(13)= 3,$ because
$n= 13, m = {{13-1}\over {v(7)}} = 4, d_7(4) = 4, d_7(13) = 7,
\varepsilon(7) = 1, n - d_7(n) = 6, m - d_7(m) = 0$
$ord_7(\omega(13)) = 4 - {6 \over 6} = 3$
\item[(2)] $ord_7 \omega(100)= 21,$ because $v(7) = 3, n = 100, m =[ {{n-1}\over {3}}]
= 33, d_7 = 100) = 4, d_7(33) = 9$ $ord_7(\omega(100)) = 33 -
{{96-24} \over {6}} = 21,$ because $log(\omega(n)) \approx {n+1
\choose 2} log(2) + n - (n + 1/2) log(n).$
\end{itemize}
\end{example}

\medskip

Let $\pi_M(x)$ be defined to be number of odd primes $p$ such that
the Mersenne order $v(p)$ is $\le x-1.$

It seems interesting to determine the asymptotic behaviour of
$\pi_M.$

One can expect some results about $\pi_M$ from the prime
factorization of $n!_M,$ because $\pi_M(n)$ is the number of
primes dividing $n!_M.$
\begin{example} $\pi_M(16) = 15$ because
the primes dividing $16!_M$ are the odd primes $\le 16$ and the
primes 17, 23, 31, 43, 73, 89, 127, 151, 257, 8191.
\end{example}
\end{section}








\begin{section}{Coefficients of the non-associative exponential}

Let $M$ be the magma with neutral element $1_M$ freely generated
by a set consisting of one element $x$. The multiplication in $M$
is a map $\cdot : M \times M \to M$ and the restriction of this
map onto $(M - \{1_M\} ) \times (M x \{1_M\})$ is injective and
its image is $M - \{1_M, x\}.$ It means that for any $t \in M -
\{1_M, x \}$ there is a unique pair $(t_1, t_2) \in M \times M, t_1
\neq 1_M \neq t_2$, such that $t = t_1 \cdot t_2.$

There is a unique homomorphism $\deg : M \rightarrow \N$ such that
$deg(1_M) = 0, deg(x) = 1.$

It is important to realize that any $t \in M, t \neq 1_M,$ gives
rise to a unique planar binary rooted tree with $deg(t)$ leaves, see
[1], section 1, p. 163. The grafting of such trees corresponds to
the multiplication map of $M$.

Let $P = \Q \{\{x \}\}$ be the $\Q$-algebra of power series with
monomials from $M$. Thus any $f \in P$ has a unique expansion $$f
= \displaystyle \sum_{t \in M} c(t) \cdot t$$ with $c(t) \in \Q.$
It will also be called the algebra of tree power series in $x$ or
the algebra of power series over $\Q$ in a nonassociative and
noncommutative variable $x$.

Define the order $ord(f)$ of $f$ to be $min \{deg(t) : c(t) \neq 0
\},$ if $f \neq 0$ and $ord(0) = + \infty.$ The following
proposition is obtained by standard calculus methods.

\begin{proposition}
\begin{enumerate}

\item[(i)] There is a unique $\Q$-linear map
${d \over {dx}} :P \rightarrow P$ such that: $${d \over {dx}} (x) =
1$$ and $${d \over {dx}} (f \cdot g) = {d \over {dx}}(f) \cdot g + f \cdot
{d \over {dx}}(g)$$  for all $f, g \in P$.

\medskip
$f' := {d \over {dx}}(f)$
is called the derivative of $f$ with respect to $x$.

\item[(ii)] Let $g \in P$ with $ord(g) \ge 1$. Then there is a unique
$\Q$-algebra homomorphism $\eta_g : P \to P$ such that $\eta_g(x)
= g$. One denotes $\eta_g(x)$ also by $f \circ g$ or $f(g)$ and
calls $\eta_g$ the substitution homomorphism induced by $g$.

\item[(iii)] There is a canonical $\Q$-algebra homomorphism $[ \ ] :
P \rightarrow \Q[[x]],$ where $\Q[[x]]$ denotes the classical
$\Q$-algebra of power series in an associative (and commutative)
variable $x$ such that $[x]=x$.

\medskip
One calls $[f]$ the classical
power series associated to $f \in  P$.
\end{enumerate}
\end{proposition}


\begin{proposition} There is a unique $f \in \Q \{\{x\}\}$ such
that $$f'(0) = 1$$ $$f(x) \cdot f (x) = f (2x)$$

Moreover $f'(x) = f(x)$ where $f'(x)$ denotes the derivative of
$f$ relative to $x$. The tree power series $f(x)$ is called the
nonassociative, noncommutative exponential series and will be
denoted in this article by $e^x$ or $exp(x)$.
\end{proposition}

\begin{proof}
Inductively we define a map $a: M \rightarrow \Q$ by putting $a(t)
:= 1,$ if $t=1_M$ or $t=x$ and otherwise $$a(t):= {a(t_1) \cdot
a(t_2) \over 2^n-2}$$ if $n = deg(t)$ and $t= t_1 \cdot t_2$ with
$t_i \in M.$

Put $f := \displaystyle \sum_{t \in M} a(t) \cdot t \in P = \Q
\{\{x\}\}.$ It is easy to check that $f(x) \cdot f(x) = f(2x):$

\begin{equation*}
\begin{split}
f(2x) &= \sum_{t \in M} a(t) \cdot 2^{deg(t)}\cdot t \\ f(x) \cdot
f(x) &= \sum_{(t_1,t_2) \in M \times M} a(t_1) a (t_2) t_1 \cdot
t_2 \\ & = 1+2a(x) +  \\& \ \ \ +  \sum_{t \in M, deg(t) \geq
2, t= (t_1,t_2)} (a(t) \cdot a(x^0)+a(x^0)a(t) + a(t_1) a(t_2))
\cdot t.
\end{split}
\end{equation*}
As $2a(t) + a(t_1) \cdot a(t_2) = 2^n a(t)$ for any $t \in M$ with
$deg(t) \geq 2, t = t_1 \cdot t_2,$ we obtain the desired
functional equation.

For the proof that $f'(x) = f(x)$ we refer to [DG].
%
%
%
%
%
%
%
\end{proof}
%
\medskip

From the proof we get that the coefficients $a(t)$ of $exp$
satisfy: $a(1) = a(x) = 1$ and $$a(t_1 \cdot t_2)= {{a(t_1)a(t_2)}
\over {2^n - 2}}$$ where $n = deg(t_1 \cdot t_2).$

\medskip
%

Let $\hat a(t) := n \cdot \omega(n) \cdot a(t).$

\begin{proposition}

Let $t_1,t_2 \in M,  n_i = deg(t_i) \ge 1$ and $t = t_1 \cdot
t_2, n = deg(t) = n_1 + n_2.$ Then $$\hat{a}(t) = ({{n-2}\over
{n_1 -1}})_M \hat{a}(t_1) \hat{a}(t_2)$$
\end{proposition}
\begin{proof}
\begin{equation*}
\begin{split}
\hat{a}(t)&= n!\cdot \omega(n) \cdot {{a(t_1)a(t_2)}\over{2^n - 2}}\\
 &= n!\cdot {{\omega(n)}\over {2M_{n-1}}}
 {{\hat{a}(t_1)}\over {n_1! \omega(n_1)}}{{\hat{a}(t_2)}\over {n_2!  \omega(n_2)}}\\
&= {{2^{n-1} F(n-1)}\over {2M_{n-1}}}
{{\hat{a}(t_1)\hat{a}(t_2)}\over{2^{n_1-1}
F(n_1-1)2^{n_2-1}F(n_2-1)}}\\
 &={{F(n-2)}\over {F(n_1-1)F(n_2-1)}} \hat{a}(t_1)
\hat{a}(t_2)\\ &={n-2 \choose n_1-1}_M \hat{a}(t_1)
\hat{a}(t_2).\\
\end{split}
\end{equation*}
\end{proof}

It follows from (3.3) that $\hat{a}(t)$ is a product of Mersenne
binomials. More precisely:

\begin{corollary}
$\hat{a}(t) = \displaystyle \prod_{a \in I(t)} {n(a) - 2 \choose
n_1(a) - 1)}_M \in \N$ where $I(t)$ is the set of inner nodes of
$t$ and $n(a)$ is the degree of $t_{\le a}$ where $t_{\le a}$ is
the tree below a which is defined to consist of all
nodes $b$ for which the simple path from $b$ to the root of
$t$ is passing through $a$. Also $n_1(a)$ is the degree of the left
factor $s_1$ of $t_{\le a}$ (such that $t_{\le a} = s_1 \cdot
s_2).$
\end{corollary}
\begin{example}
Let $T_n$ be the set of trees in M of degree $n$ defined as
follows: $T_1 = \{x\}.$

If $T_{n-1}$ is already defined, then $T_n = x \cdot T_{n-1}\cup
T_{n-1}\cdot x.$

Then $\sharp T_n = 2^{n-2}$ for $n \ge 2.$ One can show that
$\hat{a}(t) = 1$ iff $t \in T_n.$
\end{example}


\begin{corollary}
$$\displaystyle \sum_{deg(t) = n} \hat{a}(t) = \omega(n)$$
\end{corollary}
\begin{proof}
We consider the canonical algebra homomorphism

$$[\ ] : \Q \{\{ x \}\} \to \Q [[x]].$$

Then $[exp]$ is the classical series $\displaystyle
\sum^\infty_{n=0} {{x^n}\over{n!}}.$ As $[t] = x^n$ if $t \in M$
has degree $n$, it follows that $\displaystyle \sum_{deg(t) = n}
a(t) = {{1}\over{n!}}.$
\end{proof}

\bigskip

For $1 \le k \le n-1$ let

$$S_k(n) = \sum_{deg(t_1) = k \atop deg(t_2) =  n-k}\hat{a}(t_1
\cdot t_2).$$

Then $$S_k(n) = {n-2 \choose k-1}_M \displaystyle \sum_{deg(t_1) =
k \atop deg(t_2) = n-k}\hat{a}(t_1) \hat{a}(t_2)$$ $$S_k(n)={n-2
\choose k-1}_M \cdot \omega(k) \cdot \omega(n-k)$$
It follows from Corollary 3.6 that $\omega(n) = \displaystyle \sum^{n-1}_{k = 1}
S_k(n)$. Thus

\begin{corollary}
$\displaystyle \sum^{n-1}_{k=1}{n-2 \choose k-1}_M \omega(k)\cdot
\omega(n-k) = \omega(n).$
\end{corollary}

\begin{example}

$\omega(1) = \omega(2) = 1, \omega(3)= 2, \omega(4)=7, \omega(5) =
42 = 2 \cdot 3 \cdot 7$

Thus
\begin{eqnarray*}
 \omega(6) &=& \displaystyle \sum^5_{k=1} {4 \choose k-1}_M
\omega(k) \omega(6-k)\\ &=& 1 \cdot \omega(5) + M_4 \cdot 1 \cdot
\omega(4) + {{M_4 \cdot M_3} \over{M_2}} \omega(3)^2 + M_4 \cdot
\omega(4) \cdot 1 + 1 \cdot \omega(5)\\ &=& 2 \cdot 42 + 30 \cdot
7 + 5 \cdot 7 \cdot 4 = 434 = 2 \cdot 7 \cdot 31
\end{eqnarray*}

How to understand from this partition of $\omega(6)$ that 31 is a
factor of $\omega(6)?$
\end{example}

\end{section}

\bigskip\bigskip\bigskip

\pagebreak

{\bf{\large{Key Words}}}

Binary rooted trees, non-associative power series, exponential
series, $q$-binomials, Mersenne binomials

\pagebreak

{\bf{\large{Address:}}}
\bigskip
Prof. Dr. Lothar Gerritzen

Ruhr-Universit\"at Bochum

Fakult\"a f\"ur Mathematik

D 44780 Bochum

 Germany

Tel.-number: 0049-0234-32-28304

FAX: 0049-0234-32-14025

Email: Lothar.Gerritzen@rub.de
\end{document}